\documentclass[letterpaper, 10 pt, conference]{ieeeconf}

\usepackage[utf8]{inputenc}

\usepackage{amsmath,amssymb,amsfonts,amsxtra,latexsym,mathrsfs}
\usepackage{dsfont}
\usepackage{accents}
\usepackage{cancel}
\allowdisplaybreaks

\usepackage{graphicx,epsfig,fancybox}
\usepackage{fancyhdr,times}
\usepackage{array,caption}
\let\labelindent\relax
\usepackage{enumitem}
\usepackage{xcolor}
\usepackage{tikz}
\usepackage{pgfplots}
\pgfplotsset{compat=1.6}
\usepgfplotslibrary{groupplots}

\makeatletter
\newcommand{\ClearThm}[1]{%
  \expandafter\let\csname #1\endcsname\relax
  \expandafter\let\csname end#1\endcsname\relax}
\makeatother

\ClearThm{proof}
\usepackage{amsthm}

\makeatletter
\g@addto@macro\normalsize{%
  \setlength\abovedisplayskip{3pt plus 1pt minus 1pt}%
  \setlength\belowdisplayskip{3pt plus 1pt minus 1pt}%
  \setlength\abovedisplayshortskip{0pt plus 1pt}%
  \setlength\belowdisplayshortskip{3pt plus 1pt minus 1pt}%
}
\makeatother
\ClearThm{theorem}\ClearThm{lemma}\ClearThm{corollary}\ClearThm{prop}
\ClearThm{assumption}\ClearThm{definition}\ClearThm{example}
\ClearThm{exmp}\ClearThm{remark}\ClearThm{plain}\ClearThm{mydef}

\newtheorem{lemma}{Lemma}
\newtheorem{theorem}{Theorem}
\newtheorem{prop}{Proposition}[section]
\newtheorem{assumption}{Assumption}[section]

\newtheorem{corollary}{\bf{Corollary}}[section]
\newtheorem{remark}{Remark}[section]
\theoremstyle{definition}

\newtheorem{definition}{Definition}

\DeclareMathOperator{\Law}{Law}

   \def\B{{\mathcal B}}

\def\P{{\mathcal P}}      
      
\def\E{{\mathds E}}

\pgfplotsset{soldot/.style={color=blue,only marks,mark=*}}
\pgfplotsset{holdot/.style={color=blue,fill=white,only marks,mark=*}}

\fussy
\IEEEoverridecommandlockouts
\renewcommand{\baselinestretch}{0.94}

\usepackage{hyperref}

\begin{document}
\sloppy
\title{Discrete-Time Stochastic Control Beyond i.i.d. Noise
}


\author{Reyajul Hasan Leon 
     and Ali Devran Kara\thanks{The authors are with the Department of Mathematics,
     Florida State University, Tallahassee, FL, USA, email:\{rl24c,akara\}@fsu.edu}}
\maketitle
\begin{abstract}
We study discrete-time stochastic control problems where the noise driving the state dynamics is a stationary infinite-order chain rather than an i.i.d. sequence. We formulate the problem as a controlled Markov process by enlarging the state space to include the posterior distribution of the noise given the state-action history. Assuming that the sensitivity of the noise transition kernel to past noise values decays exponentially, we show that the expected total variation distance between the full-history and finite-memory noise filters decays at rate $\sqrt{N}\rho^{\sqrt{N}}$. We use this result to construct finite-memory control models and establish explicit near-optimality bounds for their policies. Finally, under suitable ergodicity conditions, we show that Q-learning converges to the optimal solution of the approximate model, so the same performance bounds apply to the learned policies.
\end{abstract}

\section{Introduction}
In stochastic control, the dynamics of the state process is generally modeled as a controlled Markov process. This assumption lets one use the very well developed theoretical toolset of Markov processes and Markov decision processes such as dynamic programming, Q-learning, and their approximation guarantees. Furthermore, this assumption in many cases is general enough to accurately model many practical systems.
However, many modern applications, e.g.  wireless scheduling, sensor networks and networked control systems,  operate under disturbance structures that are too complex to be modeled as i.i.d. For such systems, the i.i.d. assumption can lead to controllers with poor performance. In this paper, we study stochastic control problems where the driving noise is a stationary infinite chain that has a dependence structure over time. 
We consider an optimal control problem with state process $(x_n)_{n\geq 0}$ taking values in $\mathds{X}$ and a control process $(u_n)_{n\geq 0}$ taking values in $\mathds{U}$ where $\mathds{X}$ and $\mathds{U}$ are assumed to be standard Borel spaces. The dynamics of the state process is given by 
\begin{align*}
x_{n+1}=f(x_n,u_n,w_n).
\end{align*}
where $(w_n)_{n\in\mathds{Z}}$ is a noise process taking values in a standard Borel space $\mathds{W}$.

Different from the standard i.i.d. assumption on the noise process $w_n$,  we assume that the noise process $(w_n)_{n\in\mathds{Z}}$ is a stationary infinite chain. We denote the entire past up to time $n$ by 
$
 \bar{w}_n:=(w_n,w_{n-1},\dots)\in \mathds{W}^{-\mathds{N}}.
$
  We assume that the noise process is stationary with distribution $\Pi(\cdot)\in\P(\mathds{W}^{-\mathds{N}})$, and that the transitions are given by
$
w_{n+1}\sim\mathcal{K}(dw_{n+1}|\bar{w}_n)
$
so that $\Pi$ is consistent with $\mathcal{K}$ that is 
$
\Pi(w_{n+1}\in A|\bar{w}_n) = \mathcal{K}(A|\bar{w}_{n}),
$
for all $A\in\B(\mathds{W})$, $\Pi$-almost surely.
Throughout the paper, for probability measures $\mu,\nu$ on a standard Borel
space $\mathds{S}$, we define the total variation distance by
\begin{align*}
\|\mu-\nu\|_{TV}&:=\sup_{\|g\|_\infty\leq 1}\bigg| \int g\,d\mu - \int g\, d\nu\bigg|
 =|\mu-\nu|(\mathds{S}),
\end{align*}
where $|\mu-\nu|$ denotes the variation measure of the signed measure
$\mu-\nu$. In
particular $\|\mu-\nu\|_{TV}\leq2$.

We assume that the state is initialized at $n=0$ independently such that $x_0\sim \mu_0\in\P(\mathds{X})$. The more general setting where the state is initialized at $n=-\infty$ can be handled using similar arguments presented in this paper, but we adopt the present formulation for simplicity of exposition.

We assume all random variables in the system, $(x_n)_{n\geq 0}$, $(u_n)_{n\geq 0}$ and $(w_n)_{n\in\mathds{Z}}$ are defined on a common probability space $(\Omega,\mathcal{F},\mathds{P})$.
The goal is to minimize the following infinite horizon discounted cost function:
\begin{align*}
&J_\beta(\mu_0,\gamma)=\sum_{k=0}^\infty \beta^k \E_\gamma\left[c(X_k,U_k)\right]
\end{align*}
where the expectation is taken with respect to the probability measure induced by the stationary noise distribution $\Pi$, the initial state $x_0\sim \mu_0$ and the policy $\gamma$. Here, $\gamma=\{\gamma_k\}_k$ is  a policy such that each $\gamma_k$ is a control function that is measurable with respect to the $\sigma-$algebra generated by the history
$
I_k=(x_0,\dots,x_k,u_0,\dots,u_{k-1})
$
with $I_0=x_0$.

\subsection{Literature Review}
The question of how well controls designed for ideal systems perform in
real systems has been of interest over the years. A common approach has
been to design controllers that work sufficiently well for all possible
uncertain systems under some structured constraints, such as $H_\infty$
norm bounded perturbations \cite{basbern,zhou1996robust}; these methods
treat the noise as an unmodeled disturbance rather than a stochastic
process with temporal correlations. In the robust stochastic control
literature, it is commonly assumed that the true model lies within an
uncertainty set defined under various metrics
\cite{jacobson1973optimal,Iyengar2005,Ghaoui2005,kara2020robustness,%
esfahani2015}. However, these methods assume that the set of possible
models are formed by Markov models, and temporal correlations in the noise
are not covered.

In continuous time models, the study of control under nonwhite noise
processes goes back to 1960's
\cite{stratonovich1965conditional,kushner1987filtering} where they consider
systems subject to wide-bandwidth noises, and study convergence properties
as the wide-bandwidth noise converges to a white noise in an appropriate
sense.  For discrete-time models, in
\cite{hairer2009ergodic,yukselSICON2017}, {\it control-free} dynamical
models driven by stationary and temporally correlated noise are studied,
with the focus on stochastic stability and ergodicity properties. In the
controlled setting, recent works such as
\cite{MAL-097,dong2022simple,chandak2024reinforcement,borkar2025stochastic}
analyze reinforcement learning with non-Markov observations using
compressed histories  and
\cite{kara2024qlearning,kara2026reinforcement} consider control problems
with general state and cost processes without the Markov assumption.
However, these works do not model the noise process explicitly, which
prevents them from providing a complete error analysis for the learned
model.
Finally, \cite{agarwal2019online} considers a linear control problem where the noise is an arbitrary disturbance sequence rather than a stochastic process with any structure. The authors study regret minimization for a controller that uses past noise realizations and the current state, against the best linear feedback controller. This is similar to our work in relaxing standard assumptions on the noise process. However, their linear dynamics allow the noise values to be extracted exactly, whereas in our setting noise information must be inferred through a filter.

\subsection{Contributions}
In Section \ref{enlarge_section}, we introduce an enlarged state space and define a Markov control problem where the sufficient information variable is the current state $x_n$ and the posterior distribution on the noise $w_n$ given the past state and control variables. Section \ref{filter_section} is the main technical section of the paper, where we establish that if the dependence of the noise kernel on past realizations decays exponentially, then the noise filter is stable with rate $\sqrt{N}\rho^{\sqrt{N}}$ for some $\rho<1$. The analysis extends the filter stability theory for partially observed MDPs to our setting, with a new challenges:  when the noise is driven by an infinite order chain, the past state and control are correlated with future noise values at arbitrary lag, which significantly complicates the dependence structure. In Section \ref{app_section}, following a similar construction in \cite{kara2023convergence}, we define an approximate model using only a truncated history of state and action variables,  whose approximation error can be bounded using the filter stability results of Section~\ref{filter_section}. In particular, this provides explicit bounds on the performance loss incurred when a finite order model formulation, in particular a standard MDP formulation, is applied to a system with temporally correlated noise. Finally, in Section \ref{learning_section}, we show that if one uses Q-learning on the truncated state and action variables, the algorithm converges to the approximate model constructed in Section \ref{app_section}, and we provide performance loss bounds for the learned policy that decay with increasing memory under appropriate conditions on the noise transition kernel.

\section{Enlargement of the State Space}\label{enlarge_section}
In this section, we enlarge the state space and define a new control problem. We first define the following belief variable over the one-sided noise sequence 
$
\pi_n:= P(\bar{w}_n\in\cdot| I_n).
$
We also define the predictor
$
\pi^-_n:= P(\bar{w}_{n-1}\in\cdot| I_n).
$
We define the new state variable as $s_n:=(x_n,\pi_n).$ We aim to show that the new state variable forms a controlled Markov process. It is a standard result in partially observed control that the belief over the noise can be represented iteratively. In particular, we can represent  the predictor  as $\pi^-_{n+1}=F(\pi_n,x_{n+1},x_n,u_n)$ independent of the control policy $\gamma$. The filter variables can also be written iteratively as:
\begin{align}\label{F_plus}
&\pi_{n+1}(\cdot) =F^+(\pi_n,x_{n+1},x_n,u_n)\nonumber\\
&:= \int_\cdot \mathcal{K}(dw_{n+1}|\bar{w}_n)F(\pi_n,x_{n+1},x_n,u_n)(d\bar{w}_n).
\end{align}
We next note that the new state process is controlled Markov:
\begin{align*}
&P(x_{n+1},\pi_{n+1}\in\cdot | x_n,\dots,x_0,\pi_n,\dots,\pi_0,u_{n},\dots,u_0)\\
&=\int_{x_{n+1}\in\cdot}\mathds{1}_{\{F^+(\pi_n,x_{n+1},x_n,u_n)\in\cdot\}} \mathcal{T}(dx_{n+1}|x_n,u_n,\pi_n)\\
&=:\eta(ds_{n+1}|s_n,u_n)
\end{align*}
where we use the notation:
\begin{align*}
\mathcal{T}(x_{n+1}\in\cdot|x_n,u_n,\pi_n)=\int \mathds{1}_{\{f(x_n,u_n,w_n)\in\cdot\}}\pi_n(dw_n).
\end{align*}
\begin{remark}
We note that  under the usual assumption that the noise process is i.i.d., the noise belief process $\pi_n$ is simply given by the distribution of the i.i.d. variables say $\pi\in\P(\mathds{W})$. In particular, the transitions between the state variables are then given by the time-invariant kernel 
\begin{align*}
\mathcal{T}(x_{n+1}\in\cdot|x_n,u_n,\pi)=\int \mathds{1}_{\{f(x_n,u_n,w_n)\in\cdot\}}\pi(dw_n)
\end{align*}
where one drops the dependence on $\pi$ and simply denote the transition kernel by $\mathcal{T}(x_{n+1}\in\cdot|x_n,u_n)$.
\end{remark}

\subsection{Stability of the Noise Belief Updates}\label{filter_section}


In this section, we study the effect of memory truncation on the noise filter, which forms the foundation of the approximation results that are based on the finite memory feedback controls.

\begin{lemma}\label{main_lemma}
Fix $0 < M < n$. Updating the posterior on the noise 
block $w_{[n-M+1:n]}$ by conditioning on the additional 
state observations $x_{[n-M+2:n+1]}$ (and controls 
$u_{[n-M+1:n]}$) can at most double the total variation 
distance between posteriors initialized from different 
priors $\mu$ and $\nu$ on the noise:
\begin{align*}
&\mathds{E}^\mu\Big[\big\|P^\mu(w_{[n-M+1:n]}
  \mid I_{n+1}) \\
&\qquad - P^\nu(w_{[n-M+1:n]}
  \mid I_{n+1})\big\|_{TV}\Big] \\
&\leq 2\,\mathds{E}^\mu\Big[\big\|P^\mu(w_{[n-M+1:n]}
  \mid I_{n-M+1}) \\
&\qquad - P^\nu(w_{[n-M+1:n]}
  \mid I_{n-M+1})\big\|_{TV}\Big].
\end{align*}
\end{lemma}
\begin{proof}
For notational simplicity, let
\begin{align*}
 Z &:= W_{[n-M+1:n]},\\
 Y &:= \bigl(X_{[n-M+2:n+1]},U_{[n-M+1:n]}\bigr),\\
 V &:= I_{n-M+1}
    =\bigl(X_{[0:n-M+1]},U_{[0:n-M]}\bigr).
\end{align*}
Then $I_{n+1}=(V,Y)$.
Fix an admissible policy $\gamma$. Given $(Z,V)=(z,v)$, the
conditional distribution of $Y$ is the same under the priors
$\mu$ and $\nu$. Indeed, once $z$ and $v$ are fixed, the state
variables in $Y$ are generated by the state dynamics and the
controls are generated by the same policy $\gamma$. We therefore
write
\[
 P^{\mu}(dy\mid z,v)
 =
 P^{\nu}(dy\mid z,v)
 =:
 G^\gamma(dy\mid z,v).
\]

For a fixed realization $V=v$, define
\begin{align*}
 \alpha_v(dz)
 &:=
 P^\mu(Z\in dz\mid V=v),\\
 \beta_v(dz)
 &:=
 P^\nu(Z\in dz\mid V=v).
\end{align*}
The corresponding conditional joint distributions of $(Z,Y)$ are
\begin{align*}
 P_v^\mu(dz,dy)
 &:=
 \alpha_v(dz)G^\gamma(dy\mid z,v),\\
 P_v^\nu(dz,dy)
 &:=
 \beta_v(dz)G^\gamma(dy\mid z,v).
\end{align*}
Let $p_v^\mu$ and $p_v^\nu$ denote their respective $Y$-marginals.
Choose regular conditional distributions
$
 P^\mu(dz\mid y,v)
$ and 
$ P^\nu(dz\mid y,v)
$
so that
\begin{align*}
 P_v^\mu(dz,dy)
 &=
 P^\mu(dz\mid y,v)p_v^\mu(dy),\\
 P_v^\nu(dz,dy)
 &=
 P^\nu(dz\mid y,v)p_v^\nu(dy).
\end{align*}
On values of $(y,v)$ that are null under $P^\nu$, the second
conditional distribution may be extended arbitrarily in a
measurable manner.
Now we define the measure
\[
 \widetilde P_v(dz,dy)
 :=
 P^\nu(dz\mid y,v)p_v^\mu(dy).
\]
Note that $P_v^\mu$ and $\widetilde P_v$ have the same
$Y$-marginal $p_v^\mu$. Thus, by 
decomposition of total variation,
\begin{align}
 &\left\|P_v^\mu-\widetilde P_v\right\|_{TV}
 \nonumber\\=
 &\int
 \left\|
 P^\mu(Z\in\cdot\mid y,v)
 -
 P^\nu(Z\in\cdot\mid y,v)
 \right\|_{TV}
 p_v^\mu(dy).
 \label{eq:fiber-tv}
\end{align}

By the triangle inequality,
\begin{align*}
 \left\|P_v^\mu-\widetilde P_v\right\|_{TV}
 &\leq
 \left\|P_v^\mu-P_v^\nu\right\|_{TV}
 +
 \left\|P_v^\nu-\widetilde P_v\right\|_{TV}.
\end{align*}
The last two measures use the same conditional kernel
$P^\nu(dz\mid y,v)$ and differ only in their $Y$-marginals.
Therefore,
\[
 \left\|P_v^\nu-\widetilde P_v\right\|_{TV}
 =
 \left\|p_v^\nu-p_v^\mu\right\|_{TV}.
\]
Since taking a marginal cannot increase total variation,
\[
 \left\|p_v^\nu-p_v^\mu\right\|_{TV}
 \leq
 \left\|P_v^\nu-P_v^\mu\right\|_{TV}.
\]
It follows that
\begin{align}
 \left\|P_v^\mu-\widetilde P_v\right\|_{TV}
 \leq
 2\left\|P_v^\mu-P_v^\nu\right\|_{TV}.
 \label{eq:joint-bound}
\end{align}

Finally, $P_v^\mu$ and $P_v^\nu$ use the same kernel  $G^\gamma(dy\mid z,v)$ on $\alpha_v$ and $\beta_v$.Hence
$
 \left\|P_v^\mu-P_v^\nu\right\|_{TV}
 \leq
 \left\|\alpha_v-\beta_v\right\|_{TV}.$

Combining this equality with \eqref{eq:fiber-tv} and
\eqref{eq:joint-bound} gives
\begin{align*}
&\int
 \left\|
 P^\mu(Z\in\cdot\mid y,v)
 -
 P^\nu(Z\in\cdot\mid y,v)
 \right\|_{TV}
 p_v^\mu(dy)
\\
&\qquad\leq
2\left\|
 P^\mu(Z\in\cdot\mid V=v)
 -
 P^\nu(Z\in\cdot\mid V=v)
\right\|_{TV}.
\end{align*}
Integrating both sides with respect to the $P^\mu$-distribution
of $V$ and recalling that $I_{n+1}=(V,Y)$ completes the proof.
\end{proof}

\begin{definition}[Coordinate wise sensitivity]\label{def:sensitivity}
For each $j\ge 1$, define
\begin{align*}
  &\delta_j
  \;=\;
  \sup\Bigl\{
    \|\mathcal{K}(\cdot\mid\bar{w})-\mathcal{K}(\cdot\mid\bar{w}')\|_{TV}
    :\\
  &\bar{w},\bar{w}' 
    \text{ agree on all coordinates except the one $j$ steps back}
  \Bigr\}.
\end{align*}
\end{definition}

\begin{assumption}\label{decay_assmp}
There exist $C>0$ and $\rho\in(0,1)$ such that $\delta_j\leq C\rho^j$. Furthermore, $\mathcal K$ is continuous under truncation of the remote past: for every $\bar w,\bar w'$, if $\bar w^{(m)}$ agrees with $\bar w'$ in its first $m$ coordinates and with $\bar w$ thereafter, then
\[
\|\mathcal K(\cdot\mid\bar w^{(m)}) -\mathcal K(\cdot\mid\bar w')\|_{TV}\to0.
\]
Finally, assume $\frac{2C\rho}{(1-\rho)^2}<1$.
\end{assumption}

We also define the sensitivity of the noise kernel to its past beyond $k$-step:
\begin{align}\label{filt_sens}
v_k(\mathcal{K}) := \sup_{\substack{\bar{w},\bar{w}'\in\mathds{W}^{-\mathds{N}}:\\ 
w_{-j}=w'_{-j},\, j=0,\ldots,k-1}} 
\|\mathcal{K}(\cdot|\bar{w}) - \mathcal{K}(\cdot|\bar{w}')\|_{TV}
\end{align}
where $v_0(\mathcal{K})$ measures the total sensitivity of $\mathcal{K}$ to its past; the Dobrushin contraction coefficient of $\mathcal{K}$ is
$\tfrac{1}{2}v_0(\mathcal{K})$. Note that, under Assumption \ref{decay_assmp}, we have
\begin{align*}
v_k({\mathcal{K}})\leq \sum_{j=k}^\infty \delta_j \leq \frac{C \rho^k}{1-\rho}.
\end{align*}
We define:
$
w_{k:n}:=(w_n,w_{n-1},\dots,w_k).
$
When the noise is driven by an infinite-order chain, the past state and control are correlated with future noise values at arbitrary lag which makes the filter stability analysis non-trivial. We instead use finite order approximations of the noise kernel to break apart this dependence structure.  

We introduce a finite-memory approximation of the noise kernel.  Fix $M \geq 1$ and an arbitrary remote past $\bar w_{-\infty:n-M-1}$. 
For $w_{n-M:n-1} = (w_{n-M},\dots,w_{n-1})$, define
\begin{align*}
&\mathcal{K}_M(dw_n \mid w_{n-M:n-1})\\
&\quad  :=\mathcal{K}\!\left(dw_n \mid w_{n-M:n-1}, \bar w_{-\infty:n-M-1}\right).
\end{align*}
We then have 
\begin{align}\label{trunc_bound}
&\sup_{ \bar w_{-\infty:n-M-1}}\big\|\mathcal{K}_M(dw_n \mid w_{n-M:n-1}) \nonumber\\
&\quad - \mathcal{K}\!\left(dw_n \mid w_{n-M:n-1}, \bar w_{-\infty:n-M-1}\right)\big\|_{TV}\leq v_M(\mathcal{K}).
\end{align}
This then introduces a trade-off in the order parameter  $M$: larger $M$ captures longer noise memory but results in a more complex dependence between state and noise, while smaller $M$ simplifies the dependence and results in a larger approximation error.
In the following, we see that by choosing $M=\sqrt{N}$, one gets a filter stability rate that behaves like $\Delta^{\sqrt{N}}+ \sqrt{N}\rho^{\sqrt{N}}$.

Throughout this section we fix $n$ and the memory length $N$, and write
\begin{align*}
&h_k := (x_{[n-N:k]},u_{[n-N:k-1]}),\qquad n-N\leq k\leq n,
\end{align*} 
for the truncated history since time $n-N$. In particular $h_n$ is the history variable of Section \ref{app_section}, and $h_{n-N}=x_{n-N}$.
For a measure $\nu$ on $\mathds{W}^{-\mathds{N}}$, we write $P^{\nu}(\,\cdot\mid h_k)$ for the posterior obtained by placing the prior $\nu$ on $\bar{w}_{n-N}$, instead of the true conditional $P(\bar{w}_{n-N}\in\cdot\mid I_{n-N})$, and updating it using the truncated history $h_k$ and using the map $F^+$ of \eqref{F_plus}.

\begin{lemma}\label{chain_bound}

For a measure $\mu$ on $\mathds{W}^{-\mathds{N}}$, define the
$M$-step joint distributions
\begin{align*}
  &\Gamma_\mu(dw_{[n-M:n-1]})\\
  &:= \int \prod_{k=n-M}^{n-1}
       \mathcal{K}(dw_k\mid\bar{w}_{k-1})\,\mu(d\bar{w}_{n-M-1}),\\
  &\Gamma_\mu^{(M)}(dw_{[n-M:n-1]})\\
  &:= \int \prod_{k=n-M}^{n-1}
       \mathcal{K}_M(dw_k\mid w_{k-M:k-1})\,\mu(d\bar{w}_{n-M-1}).
\end{align*}
For any $\|g\|_\infty \leq 1$,
\[
  \left|\int g\,d\Gamma_\mu
  - \int g\,d\Gamma_\mu^{(M)}\right|
  \leq M\, v_M.
\]
\end{lemma}
\begin{proof}
The result is a direct implication of the bound in \eqref{trunc_bound}.
\end{proof}

\begin{lemma}\label{dob_lemma}
Let $\mu,\nu$ be probability measures on $\mathds{W}^{-\mathds{N}}$.
Under Assumption~\ref{decay_assmp}, the $M$-step truncated
kernel $\mathcal{K}_M$ contracts in total variation:
\begin{align*}
  &\bigl\|\Gamma_\mu^{(M)} - \Gamma_\nu^{(M)}\bigr\|_{TV}
  \leq C\sum_{k=1}^M k\rho^k \,\|\mu-\nu\|_{TV}\\
   &\leq \frac{C\rho}{(1-\rho)^2}\,\|\mu-\nu\|_{TV}.
\end{align*}
\end{lemma}
\begin{proof}
The proof follows by a telescoping argument, adding and subtracting terms that differ in a single coordinate, and then applying Definition~\ref{def:sensitivity} and Assumption~\ref{decay_assmp}.
\end{proof}

We then have the following filter stability result, stated for an arbitrary prior assigned at time $n-N$ as a function of the state $x_{n-N}$.

\begin{prop}\label{filt_stab_prop}
Fix $N\geq 1$ and $1\leq M\leq N$. Let Assumption \ref{decay_assmp} hold, and let $\Delta:=\frac{2C\rho}{(1-\rho)^2}<1$. For $\bar{w}_0\sim \Pi$ and for $x_0$ initialized independently, we have for every collection $\{\varrho_x\}_{x\in\mathds{X}} \subset \P(\mathds{W}^{-\mathds{N}})$ that
\begin{align*}
&\sup_\gamma\E^{\Pi,\gamma}\bigg[\big\|P(w_n\in\cdot\mid I_n)  -P^{\varrho_{x_{n-N}}}(w_n\in\cdot\mid h_n) \big\|_{TV} \bigg]\\
&\leq  \frac{\Delta^{\lfloor N/M \rfloor + 1} }{2}+ \frac{(M+2)v_M}{1-\Delta}.
\end{align*}
In particular, choosing $M=\lfloor\sqrt{N}\rfloor$, the bound decays as $\frac{\Delta^{\left\lfloor \sqrt{N}\right\rfloor+1 }}{2} + O(\sqrt{N}\rho^{\sqrt{N}})$.
\end{prop}

\begin{proof}
For a measure $\mu$ on $\mathds{W}^{-\mathds{N}}$, recall the $M$-step joint distributions generated by the exact noise kernel $\mathcal{K}$ and the truncated approximate kernel $\mathcal{K}_M$ denoted by $\Gamma_\mu(dw_{[n-M:n-1]})$ and $\Gamma_\mu^{(M)}$ defined in Lemma \ref{chain_bound}.
We also introduce the following notation for simplicity
\begin{align*}
  &D_M(a,b) := \sup_\gamma\mathds{E}^{\Pi,\gamma}\Big[\big\|P(w_{[a:a+M-1]}\in\cdot \mid I_b)\\
  &\qquad\qquad - P^{\varrho_{x_{n-N}}}(w_{[a:a+M-1]}\in\cdot \mid h_b)\big\|_{TV}\Big].
\end{align*}
We  note that under Assumption \ref{decay_assmp},
\begin{align}\label{v_zero}
v_0(\mathcal{K}) \leq \sum_{j\geq 1}\delta_j \leq \frac{C\rho}{1-\rho}
 = \frac{(1-\rho)\Delta}{2}\leq \frac{\Delta}{2}.
\end{align}
Replacing $\mathcal{K}$ by $\mathcal{K}_M$ in the transition to $w_n$ from
$w_{[n-M:n-1]}$ gives an error of at most $v_M$ on each side by \eqref{trunc_bound}, and the remaining difference is bounded by the Dobrushin coefficient $\frac{1}{2}v_0(\mathcal{K})$ of $\mathcal{K}_M$. Hence we can write
\begin{align}\label{init_bound}
  &\big\|P(w_n \in \cdot \mid I_n) - P^{\varrho_{x_{n-N}}}(w_n \in \cdot \mid h_n)\big\|_{TV} \nonumber\\
  &\leq 2v_M +\frac{1}{2} v_0(\mathcal{K})\, \big\|P(w_{[n-M:n-1]} \in \cdot \mid I_n) \nonumber\\
  &\qquad - P^{\varrho_{x_{n-N}}}(w_{[n-M:n-1]} \in \cdot \mid h_n)\big\|_{TV}.
\end{align}

Using Lemma~\ref{main_lemma} with $n$ replaced by $n-1$, and with the  priors  taken to be the true conditional and $\varrho_{x_{n-N}}$ respectively, we can remove the conditioning on $x_{[n-M+1:n]}$ and $u_{[n-M:n-1]}$ to write:
\begin{align}\label{bound_lem1}
  D_M(n-M, n) \leq 2\, D_M(n-M, n-M).
\end{align}
Note that the expectation in Lemma~\ref{main_lemma} is taken under the first prior. In the above, it corresponds to the true conditional, so that the expectation operator is $\E^{\Pi,\gamma}$ as needed.

Next, we compare the exact and truncated $M$-step joint distributions. Writing $\mu = P(\bar{w}_{n-M-1} \in \cdot \mid I_{n-M})$ for the posterior under the
full history, Lemma~\ref{chain_bound} gives $\|\Gamma_\mu - \Gamma_\mu^{(M)}\|_{TV} \leq Mv_M$. Writing $\nu = P^{\varrho_{x_{n-N}}}(\bar{w}_{n-M-1} \in \cdot \mid h_{n-M})$ for the posterior under the truncated history and the prior $\varrho_{x_{n-N}}$, the same bound gives $\|\Gamma_\nu - \Gamma_\nu^{(M)}\|_{TV} \leq Mv_M$. Once both sides use $\mathcal{K}_M$, Lemma~\ref{dob_lemma} can be used to get:
\begin{align*}
  \big\|\Gamma_\mu^{(M)} - \Gamma_\nu^{(M)}\big\|_{TV}\leq \frac{C\rho}{(1-\rho)^2}\,\|\mu - \nu\|_{TV}.
\end{align*}
Since $\Gamma^{(M)}_\mu$ and $\Gamma^{(M)}_\nu$ depend on $\mu$ and $\nu$ only through their marginals on $w_{[n-2M:n-M-1]}$,
so the right hand side may be taken with $\|\mu-\nu\|_{TV}$ replaced by the total variation distance between those marginals, whose expectation is
$D_M(n-2M,n-M)$. Combining these bounds:
\begin{align*}
 &D_M(n-M, n-M) \\
 &\leq 2Mv_M + \frac{C\rho}{(1-\rho)^2}\,D_M(n-2M, n-M).
\end{align*}
Together with \eqref{bound_lem1}:
\begin{align}\label{one_block}
  D_M(n-M, n) \leq 4Mv_M + \Delta\, D_M(n-2M, n-M).
\end{align}

Repeating \eqref{one_block}  $\lfloor N/M \rfloor$ times, the second
argument reaches $b=n-N$, where $h_{n-N}=x_{n-N}$ and the two measures being
compared are the true conditional $P(\bar{w}_{n-N}\in\cdot\mid I_{n-N})$ and
the prior $\varrho_{x_{n-N}}$. Using the trivial bound of $2$ on that TV
difference we get
\begin{align*}
  D_M(n-M, n) \leq 2\Delta^{\lfloor N/M \rfloor}  + \frac{4Mv_M}{1-\Delta}.
\end{align*}
Substituting back into \eqref{init_bound}, and using \eqref{v_zero} :
\begin{align*}
 &\sup_\gamma \mathds{E}^{\Pi,\gamma}\bigg[\big\|P(w_n \in \cdot \mid I_n)
  - P^{\varrho_{x_{n-N}}}(w_n \in \cdot \mid h_n)\big\|_{TV}\bigg] \\
 &\leq   \frac{1}{2}\Delta^{\lfloor N/M \rfloor + 1} + \frac{(M+2)v_M}{1-\Delta}.
\end{align*}
\end{proof}

We define the loss constant for the future results:
\begin{align}\label{loss_bnd}
\mathcal{E}_N := \tfrac{1}{2}\Delta^{\lfloor\sqrt{N}\rfloor+1}+ \frac{\bigl(\lfloor\sqrt{N}\rfloor+2\bigr)C\rho^{\lfloor\sqrt{N}\rfloor}}  {(1-\rho)(1-\Delta)}.
\end{align}

\begin{corollary}\label{filt_stab_cor}
Under the hypotheses of Proposition \ref{filt_stab_prop}, taking
$\varrho_x = P(\bar{w}_{n-N}\in\cdot\mid x_{n-N}=x)$ gives 
\begin{align*}
&\sup_\gamma\E^{\Pi,\gamma}\Big[\big\|P(w_n\in\cdot\mid I_n)
 -P(w_n\in\cdot\mid h_n) \big\|_{TV} \Big] \leq \mathcal{E}_N.
\end{align*}
\end{corollary}

\section{Near Optimality of Lower Order Approximations}\label{app_section}
Recall the controlled Markov model we introduced in Section \ref{enlarge_section} with the state $s_n=(x_n,\pi_n)$ where $\pi_n=P(\bar{w}_n\in\cdot| I_n)$ denotes the filter variable for the one-sided noise sequence. 

In this section, following a similar approach to \cite{kara2023convergence}, we will introduce an alternative state representation without loss of optimality, that is more suitable for  finite-memory based approximations.  For a given memory length $N<\infty$, we let the new state be the pair:
\begin{align*}
&h_n  = (x_{[n-N:n]},u_{[n-N:n-1]})\\
&\pi_{n-N}= P(\bar{w}_{n-N}\in\cdot\mid I_{n-N})
\end{align*}
where $h_n$ is the history variable introduced in Section \ref{filter_section}.
For the transitions from the pair $h_n, \pi_{n-N}$ under the control action $u_n$, to $h_{n+1},\pi_{n-N+1}$, note that $\pi_{n-N+1}$ is determined by the map $F^+(\pi_{n-N},x_{n-N+1},x_{n-N},u_{n-N})$ defined in (\ref{F_plus}), hence the only randomness is on the evolution of $x_{n+1}$:
\begin{align}\label{x_eval}
&P(x_{n+1}\in\cdot|h_n,\pi_{n-N},u_n)\nonumber\\
&\quad=\int \mathds{1}_{\{f(x_n,u_n,w_n)\in\cdot \}} P(dw_n| h_n,\pi_{n-N}).
\end{align}
For the following we define 
$
s_n=(h_n,\pi_{n-N})
$
to be the new state variable. We denote the conditional expectation of a function $g$ of the state variable $s$ from some state $s_n=(h_n,\pi_{n-N})$ and control $u_n$ by
\begin{align*}
\int_{x_{n+1}\in\mathds{X}} g(s_{n+1})P(dx_{n+1}|h_n,\pi_{n-N},u_n)
\end{align*}
since the only randomness is in the transition of $x_{n+1}$ and all other components of $s_{n+1}$ evolve deterministically.

\subsection{Approximate Model Construction}
We now define the approximate model, bound the single-step approximation
error, and state the main performance guarantee. For the approximate model, we fix a collection  $\varrho=\{\varrho_x\}_{x\in\mathds{X}}\subset \mathcal{P}(\mathds{W}^{-\mathds{N}})$, where for each $x\in\mathds{X}$, $\varrho_x\in\mathcal{P}(\mathds{W}^{-\mathds{N}})$  assigns a prior on the noise past. We note that $\varrho$ is a  design choice, and the results of this section hold for every such choice. We define a quantization map $\phi$  by
\begin{align*}
\phi(s_n) = \phi(h_n, \pi_{n-N}) = (h_n, \varrho_{x_{n-N}}),
\end{align*}
where $x_{n-N}$ is the oldest state coordinate of the history variable $h_n$. The map $\phi$ effectively discards the dependence on state and control variables more than $N$ steps in the past, keeping only $x_{n-N}$. In particular, choosing $N=0$ reduces the model to a standard MDP. In Section \ref{learning_section}, we will observe that when the model is learned from an ergodic path of data, the learned transition kernel corresponds to the particular choice in which $\varrho_x$ is the conditional distribution of $\bar{w}_{n-N}$ given $x_{n-N}=x$ under the stationary distribution of the pair $(\bar{w}_n,x_n)$.
For a test function $g$ of the original state variable $s_n = (h_n, \pi_{n-N})$, we define the approximate transition kernel $\hat{P}$ by
\begin{align*}
\hat{P}(x_{n+1} \mid s_n, u_n) = P(x_{n+1} \mid \phi(s_n), u_n).
\end{align*}
That is, the approximate model $\hat{P}$ replaces the true conditional on $s_n$ with the conditional on the quantized state $\phi(s_n)$. In particular, using the notation introduced before Proposition \ref{filt_stab_prop}, 
\begin{align}\label{approx_kernel}
&\hat{P}(x_{n+1}\in\cdot \mid s_n,u_n)\nonumber\\
&\quad=\int \mathds{1}_{\{f(x_n,u_n,w_n)\in\cdot\}}\,
   P^{\varrho_{x_{n-N}}}(dw_n \mid h_n).
\end{align}
The right hand side of \eqref{approx_kernel} depends on $s_n$  through $h_n$, and $h_{n+1}$ is determined by $(h_n,u_n,x_{n+1})$. Hence, the approximate model is Markov on the finite history variables $h_n$. We denote the value function of the approximate model by $J_\beta^N$ and its optimal policy by $\gamma^N$.

We now compare the original and approximate dynamics. For any $\|g\|_\infty \leq 1$, using the dynamics of $x_n$ given in \eqref{x_eval}, 
\begin{align}\label{noise_TV_bound}
&\bigg|\int_{x_{n+1}} g(s_{n+1})\,P(dx_{n+1} \mid \phi(s_n), u_n) \nonumber\\
&\quad - \int_{x_{n+1}} g(s_{n+1})\,P(dx_{n+1} \mid s_n, u_n)\bigg| \nonumber\\
&\leq \bigl\| P^{\varrho_{x_{n-N}}}(dw_n \mid h_n) - P(dw_n \mid I_n) \bigr\|_{TV}.
\end{align}
For the last term, we have the $TV$-distance of two posteriors on the noise variable $w_n$,  the first one starts from the prior $\varrho_{x_{n-N}}$ updated using the truncated
history and the second one is the true filter given the full history. Proposition \ref{filt_stab_prop} can be used to bound this in expectation, uniformly over the
choice of $\varrho$.

\begin{theorem}\label{approx_model_thm}
Suppose the system is controlled under an arbitrary policy $\gamma$ for the time steps $n = 0, \dots, N-1$, and consider the starting state $s_N = (h_N, \pi_0)$ with
\begin{align*}
&h_N = (x_{[0:N]},u_{[0:N-1]}), \\
&\pi_0 = P(\bar{w}_0\in\cdot \mid x_0) = \Pi.
\end{align*}
Let $J_\beta^N$ denote the optimal value function of the approximate model and let $\gamma^N$ denote its corresponding optimal policy. 
Under Assumption \ref{decay_assmp}, the following bounds hold:
\begin{align*}
&\mathds{E}\Big[\big|J_\beta^*(s_N) - J_\beta^N(s_N)\big|\Big]\leq \frac{\|c\|_\infty\, \beta}{(1-\beta)^2}
\mathcal{E}_N, \\[6pt]
&\mathds{E}\Big[\big|J_\beta(s_N, \gamma^N) - J_\beta^*(s_N)\big|\Big]\leq \frac{2\|c\|_\infty\, \beta}{(1-\beta)^2}
\mathcal{E}_N,
\end{align*}
where $\mathcal{E}_N$ is defined in \eqref{loss_bnd} and the expectation is taken over the random realizations of $s_N = (h_N,\pi_0)$, with the controls during $n = 0, \dots, N-1$ chosen according to $\gamma$.
\end{theorem}
\begin{proof}
We define 
\begin{align*}
&L(s_n):= \sup_u\big\|P(dx_{n+1}|s_n,u) - P(dx_{n+1}|\phi(s_n),u) \big\|_{TV}.
\end{align*}
Using the Bellman equations for $J_\beta^*$ and $J^N_\beta$ we then have that 
\begin{align*}
&\big|J_\beta^*(s_N) - J^N_\beta(s_N) \big| \\
&\leq \sup_u\bigg| \beta \int J_\beta^*(s_{N+1})P(dx_{N+1}|s_N,u) \\
&\quad-  \beta \int J_\beta^N(\phi(s_{N+1}))P(dx_{N+1}|\phi(s_N),u) \bigg|\\
&\leq \beta \|J_\beta^N\|_\infty L(s_N) \\
&\quad + \sup_\gamma\beta \mathds{E}\bigg[\big|J_\beta^*(s_{N+1}) - J^N_\beta(s_{N+1}) \big| \mid s_N,\gamma(s_N)\bigg].
\end{align*}
Iterating this inequality, at each step the pointwise supremum over control functions is bounded by the supremum over policies, so that taking the expectation we write:
\begin{align*}
&\mathds{E}\bigg[\big|J_\beta^*(s_N) - J^N_\beta(s_N) \big|\bigg] \\
& \leq \beta\|J_\beta^N\|_\infty   \sup_\gamma\mathds{E}^\gamma\bigg[\sum_{n=0}^\infty \beta^nL(s_{N+n})\bigg]
\end{align*}
where the expectation  is with respect to the true dynamics without memory truncation. Using the bound in \eqref{noise_TV_bound} and using Proposition \ref{filt_stab_prop}, we can conclude the first part of the proof.
For the second part, we write
\begin{align*}
&|J_\beta(s_N,\gamma^N) - J_\beta^*(s_N) |\\
&\quad\leq |J_\beta(s_N,\gamma^N) - J_\beta^N(s_N)| +  |J_\beta^N(s_N)- J_\beta^*(s_N) |
\end{align*}
the second term is bounded in the first part of the proof, and the first part can be bounded using identical arguments, since both terms are controlled by the same policy $\gamma^N$.
\end{proof}

\begin{remark}\label{rmk:robustness}
Theorem~\ref{approx_model_thm} can be seen as a robustness bound. It gives an upper bound on the performance loss when a finite-memory approximation of the noise  kernel is used instead of the true infinite order model. In particular, choosing $N=0$ gives a standard MDP with transitions
\[
  P(x_{n+1} \mid x_n, u_n) = \int \mathds{1}_{\{x_{n+1} = f(x_n,u_n,w_n)\}}\,\varrho_{x_n}(dw_n),
\]
whose noise law is allowed to depend on the current state. Taking in addition $ \varrho_x\equiv\Pi$ for every $x$, where $\Pi$ denotes the marginal of the stationary noise distribution, recovers the standard MDP built from the unconditional noise law. Applying Proposition \ref{filt_stab_prop} with $N=0$
and any $M\geq 1$ gives the bound
$\frac{\Delta}{2}+\frac{(M+2)C\rho^{M}}{(1-\rho)(1-\Delta)}$, and letting
$M\to\infty$ the performance loss is bounded by
\[
\frac{\|c\|_\infty\,\beta}{(1-\beta)^2}\,\frac{\Delta}{2},
\qquad \Delta=\frac{2C\rho}{(1-\rho)^2}.
\]
Thus, the MDP approximation performs well when $\rho$, the dependence decay
rate of the noise kernel, is close to zero.
\end{remark}

\subsection*{ Discussion on the approximation approach. }
An alternative approach would be to first approximate the noise process by a finite-order model governed by the truncated kernel $\mathcal{K}_N$, and then solve the resulting control problem. However, even after such a reduction, the noise remains past dependent (as long as $N\geq 1$), and full optimality still  requires the full history of observable state variables. Finite-memory approximations on the state history are therefore still needed.  Furthermore, when the noise dynamics as well as the state dynamics are unknown and one learns the model directly from state  observations, the learned model naturally  corresponds to a finite-memory approximation of  the infinite-order noise, since finite order  noise reduction requires knowledge of the noise kernel. For these reasons, we construct finite memory approximations directly on the observable history. We note that finite-order approximations of the noise are used  indirectly in our analysis, as a technical tool in the filter stability proof (Proposition~\ref{filt_stab_prop}).

\section{Learning}\label{learning_section}
We assume in this section that $\mathds{X}$ and $\mathds{U}$ are finite. Consider the following Q-learning iterations 
\begin{align}\label{iterateAlgM}
&Q_{t+1}(h,u)=\left(1-\alpha_t(h,u)\right)Q_t(h,u) \nonumber\\
&\quad + \alpha_t(h,u)\left(c(x_t,u_t) + \beta V_t(h_{t+1}) \right)
\end{align}
for all $(h,u)\in \mathds{X}^{N+1}\times\mathds{U}^N\times\mathds{U}$ where $V_t(h)=\min_u Q_t(h,u)$.
We assume the data is collected under a stationary randomized Markov exploration policy $\gamma_{\mathrm e}(\cdot\mid x)$ with $\gamma_{\mathrm e}(u\mid x)>0$ for all $(x,u)$.
\begin{assumption}\label{erg_assmp}
The joint process $(h_{t+1},h_t,u_t,\bar{w}_t)_{t\geq 0}$ is asymptotically ergodic in the sense that {for the given initialization $h_0$}, for any measurable bounded function $f$, we have that with probability one,
\begin{align*}
\frac{1}{T}\sum_{t=0}^{T-1} f(h_{t+1},h_t,u_t,\bar{w}_t)\to \int f\,d\pi
\end{align*}
for some probability measure $\pi$. Furthermore, the marginal of $\pi$ on $(h,u)$ is strictly positive, i.e., $\pi(h,u)>0$ for all $h\in \mathds{X}^{N+1}\times\mathds{U}^N$ and $u\in \mathds{U}$.
\end{assumption}

By \cite[Theorem 2.1]{kara2024qlearning}, the iterations \eqref{iterateAlgM} converge to the optimal $Q$-function of the finite MDP with kernel
$\pi(x_{t+1}=x'\mid h_t=h,u_t=u)$. We claim this is the approximate kernel \eqref{approx_kernel} with
$
  \varrho_x := \pi\big(\bar{w}_{t-N}\in\cdot \mid x_{t-N}=x\big),
$
which is well defined since $\mathds{X}$ is finite and $\pi(x_{t-N}=x)>0$.

The relations $x_{k+1}=f(x_k,u_k,w_k)$, $\Law(w_{k+1}\mid\bar{w}_k)=\mathcal{K}(\cdot\mid\bar{w}_k)$ and $\Law(u_k\mid x_k)=\gamma_{\mathrm e}(\cdot\mid x_k)$ hold along the trajectory, and thus they hold under the regular conditionals under $\pi$. Therefore, under $\pi$, the variables over the finite-window evolve according to the true dynamics, and the conditional distribution of $\bar{w}_n$ given $h_n$ is given by the Bayesian recursion of Section \ref{enlarge_section}, started from the conditional law of $\bar{w}_{n-N}$ given the information available at time $n-N$ in the window, which is only $x_{n-N}$. Hence
$
  \Law_\pi\big(\bar{w}_n\mid h_n\big)
  = P^{\varrho_{x_{n-N}}}\big(\,\cdot\mid h_n\big),
$
and since the exploration policy depends only on the current state, so
$\pi(x_{n+1}=x'\mid h_n,u_n)=\hat{P}(x'\mid h_n,u_n)$ by
\eqref{approx_kernel}.

\begin{corollary}
Under Assumption \ref{erg_assmp}  if $\sum_t\alpha_t(h,u)=\infty$ and $\sum_t \alpha^2_t(h,u)<\infty$ a.s., the iterations in \eqref{iterateAlgM}converge almost surely to some $Q(h,u)$, where $Q(h,u)$ is the optimal Q-value function of the approximate model constructed in Section \ref{app_section}.
Furthermore, suppose the system is controlled under an arbitrary policy $\gamma$ for the time steps $n = 0, \dots, N-1$, and consider the starting state $s_N = (h_N, \pi_0)$ with
\begin{align*}
&h_N = (x_{[0:N]},u_{[0:N-1]}),\quad \pi_0 = P(\bar{w}_0\in\cdot\mid x_0)=\Pi.
\end{align*}
Under Assumption \ref{decay_assmp}, denoting the learned policy by $\gamma^N$, the bound of Theorem \ref{approx_model_thm} then applies to $\gamma^N$.
\end{corollary}

\section{Conclusion}
We studied  control problems driven by exogenous noise processes  given by an infinite order chain. The main question addressed in this paper is whether, and under what conditions, finite-order model approximations, and in particular controllers designed under the Markov assumption, are robust when they are used in environments  where the noise has temporal correlations. We established that if the dependence of the noise kernel on its past decays exponentially, then one can provide approximation and performance loss bounds that decay with the approximation order.

\bibliographystyle{IEEEtran}

\bibliography{SerdarBibliography,AliBibliography,references_acc,references}

\begin{thebibliography}{10}
\providecommand{\url}[1]{#1}
\csname url@samestyle\endcsname
\providecommand{\newblock}{\relax}
\providecommand{\bibinfo}[2]{#2}
\providecommand{\BIBentrySTDinterwordspacing}{\spaceskip=0pt\relax}
\providecommand{\BIBentryALTinterwordstretchfactor}{4}
\providecommand{\BIBentryALTinterwordspacing}{\spaceskip=\fontdimen2\font plus
\BIBentryALTinterwordstretchfactor\fontdimen3\font minus
  \fontdimen4\font\relax}
\providecommand{\BIBforeignlanguage}[2]{{%
\expandafter\ifx\csname l@#1\endcsname\relax
\typeout{** WARNING: IEEEtran.bst: No hyphenation pattern has been}%
\typeout{** loaded for the language `#1'. Using the pattern for}%
\typeout{** the default language instead.}%
\else
\language=\csname l@#1\endcsname
\fi
#2}}
\providecommand{\BIBdecl}{\relax}
\BIBdecl

\bibitem{basbern}
T.~Ba\c{s}ar and P.~Bernhard, \emph{H-infinity Optimal Control and Related
  Minimax Design Problems: A Dynamic Game Approach}.\hskip 1em plus 0.5em minus
  0.4em\relax Boston, MA: Birkh\"auser, 1995.

\bibitem{zhou1996robust}
K.~Zhou, J.~C. Doyle, and K.~Glover, \emph{Robust and optimal control}.\hskip
  1em plus 0.5em minus 0.4em\relax Prentice-Hall, 1996, vol.~40.

\bibitem{jacobson1973optimal}
D.~Jacobson, ``Optimal stochastic linear systems with exponential performance
  criteria and their relation to deterministic differential games,'' \emph{IEEE
  Transactions on Automatic control}, vol.~18, no.~2, pp. 124--131, 1973.

\bibitem{Iyengar2005}
G.~N. Iyengar, ``Robust dynamic programming,'' \emph{Mathematics of Operations
  Research}, vol.~30, no.~2, pp. 257--280, 2005.

\bibitem{Ghaoui2005}
A.~Nilim and L.~E. Ghaoui, ``Robust control of {M}arkov decision processes with
  uncertain transition matrices,'' \emph{Operations Research}, vol.~53, no.~5,
  pp. 780--798, 2005.

\bibitem{kara2020robustness}
A.~Kara and S.~Y\"uksel, ``Robustness to incorrect system models in stochastic
  control,'' \emph{SIAM Journal on Control and Optimization}, vol.~58, no.~2,
  pp. 1144--1182, 2020.

\bibitem{esfahani2015}
P.~M. Esfahani and D.~Kuhn, ``Data-driven distributionally robust optimization
  using the {W}asserstein metric: Performance guarantees and tractable
  reformulations,'' \emph{Mathematical Programming}, pp. 1--52, 2017.

\bibitem{stratonovich1965conditional}
R.~L. Stratonovich, ``Conditional markov processes,'' in \emph{Non-linear
  transformations of stochastic processes}.\hskip 1em plus 0.5em minus
  0.4em\relax Elsevier, 1965, pp. 427--453.

\bibitem{kushner1987filtering}
H.~Kushner and W.~Runggaldier, ``Filtering and control for wide bandwidth noise
  driven systems,'' \emph{IEEE Transactions on Automatic Control}, vol.~32,
  no.~2, pp. 123--133, 1987.

\bibitem{hairer2009ergodic}
M.~Hairer, ``Ergodic properties of a class of non-markovian processes,''
  \emph{Trends in stochastic analysis}, vol. 353, pp. 65--98, 2009.

\bibitem{yukselSICON2017}
S.~Y{\"u}ksel, ``On stochastic stability of a class of non-{M}arkovian
  processes and applications in quantization,'' \emph{SIAM J. on Control and
  Optimization}, vol.~55, pp. 1241--1260, 2017.

\bibitem{MAL-097}
\BIBentryALTinterwordspacing
X.~Lu, B.~V. Roy, V.~Dwaracherla, M.~Ibrahimi, I.~Osband, and Z.~Wen,
  ``Reinforcement learning, bit by bit,'' \emph{Foundations and Trends® in
  Machine Learning}, vol.~16, no.~6, pp. 733--865, 2023. [Online]. Available:
  \url{http://dx.doi.org/10.1561/2200000097}
\BIBentrySTDinterwordspacing

\bibitem{dong2022simple}
S.~Dong, B.~van Roy, and Z.~Zhou, ``Simple agent, complex environment:
  Efficient reinforcement learning with agent states,'' \emph{The Journal of
  Machine Learning Research}, vol.~23, no.~1, pp. 11\,627--11\,680, 2022.

\bibitem{chandak2024reinforcement}
S.~Chandak, P.~Shah, V.~S. Borkar, and P.~Dodhia, ``Reinforcement learning in
  non-{M}arkovian environments,'' \emph{Systems \& Control Letters}, vol. 185,
  p. 105751, 2024.

\bibitem{borkar2025stochastic}
V.~S. Borkar, ``Stochastic approximation in non-markovian environments,''
  \emph{Systems \& Control Letters}, vol. 205, p. 106250, 2025.

\bibitem{kara2024qlearning}
\BIBentryALTinterwordspacing
A.~D. Kara and S.~Yuksel, ``Q-learning for stochastic control under general
  information structures and non-markovian environments,'' \emph{Transactions
  on Machine Learning Research}, 2024, featured Certification. [Online].
  Available: \url{https://openreview.net/forum?id=1Yp6xpTV55}
\BIBentrySTDinterwordspacing

\bibitem{kara2026reinforcement}
A.~D. Kara, ``Reinforcement learning with function approximation for non-markov
  processes,'' \emph{arXiv preprint arXiv:2601.00151}, 2026.

\bibitem{agarwal2019online}
N.~Agarwal, B.~Bullins, E.~Hazan, S.~M. Kakade, and K.~Singh, ``Online control
  with adversarial disturbances,'' in \emph{Proc. 36th Int. Conf. Mach. Learn.
  (ICML)}, 2019, pp. 111--119.

\bibitem{kara2023convergence}
A.~D. Kara and S.~Y{\"u}ksel, ``Convergence of finite memory {Q} learning for
  {POMDP}s and near optimality of learned policies under filter stability,''
  \emph{Mathematics of Operations Research}, vol.~48, no.~4, pp. 2066--2093,
  2023.

\end{thebibliography}

\end{document}